\documentclass{elsart}

\usepackage{amsfonts}
\usepackage{amssymb}
\usepackage{latexsym}

\newcommand{\R}{\mathbb{R}}
\newcommand{\Z}{\mathbb{Z}}

\begin{document}

\begin{frontmatter}

\title{Characteristic Classes of Bad Orbifold Vector Bundles
\thanksref{fde}}
\thanks[fde]{partially supported by a Rhodes College Faculty Development Endowment Grant}

\author{Christopher Seaton}
\address{Department of Mathematics and Computer Science,
Rhodes College, 2000 N. Parkway,
Memphis, TN 38112, USA.
Tel: (901)843-3721,
Fax: (901)843-3050,
e-mail: seatonc@rhodes.edu }

\begin{abstract}

We show that every bad orbifold vector bundle can be realized as the restriction
of a good orbifold vector bundle to a suborbifold of the base space.  We give an
explicit construction of this result in which the Chen-Ruan orbifold cohomology
of the two base spaces are isomorphic (as additive groups).  This construction is used
to indicate an extension of the Chern-Weil construction of characteristic classes
to bad orbifold vector bundles.  In particular, we apply this construction to the
orbifold Euler class and demonstrate that it acts as an obstruction to the existence
of nonvanishing sections.

\end{abstract}


\end{frontmatter}

% XXXXXXXXXXXXXXXXXXXXXXXXXXXXXXXXXXXXXXXXXXXXXXXXXXXXXXXXXXXXXX
% XXXXXXXXXX            SECTION: Introduction
% XXXXXXXXXXXXXXXXXXXXXXXXXXXXXXXXXXXXXXXXXXXXXXXXXXXXXXXXXXXXXX

\section{Introduction}

When orbifolds were originally introduced by Satake (under the name
$V$-manifold; see \cite{satake1} and \cite{satake2}), the definition
given coincides with the modern definition of a reduced (codimension
2) orbifold. Recall that an orbifold $Q$ is a Hausdorff topological
space locally modeled on $\R^n / G$ where $G$ is a finite group of
automorphisms; $Q$ is {\bf reduced} if each local group $G$ acts
effectively and {\bf unreduced} otherwise.  Since each point of an
unreduced orbifold is a fixed point for a nontrivial group element,
the orbifold is composed entirely of singular points.  Such
orbifolds appear, for instance, as suborbifolds of reduced orbifolds
consisting entirely of singular points.

Many of the techniques used to study the differential geometry of vector bundles over smooth
manifolds extend readily to the case of orbifold vector bundles over smooth orbifolds provided
that the orbifolds in question are reduced  (see \cite[Section 4.3]{ruangwt}).  In this case,
every orbifold vector bundle is a good orbifold vector bundle.  An orbifold vector bundle $E$
over the orbifold $Q$ is {\bf good} if it is covered by charts of
the form $\{ V \times \R^k , G, \tilde{\pi} \}$ such that $\{ V,
\R^k, \pi \}$ is an orbifold chart over $Q$, and the kernel of the
$G$-action on $V \times \R^k$ coincides with the kernel of the
$G$-action on $V$.  If $E$ is a good vector bundle, then its geometry is identical
to that of a vector bundle over a reduced orbifold.  In fact, as noted by Chen and Ruan
\cite{ruangwt}, the associated reduced
orbifold $E_{red}$ is an orbifold vector bundle over $Q_{red}$.
Therefore, the fibers of $E$ over an open dense subset of $Q$ are vector spaces.

In the case that $E$ is not good, it is called a {\bf bad} orbifold
vector bundle. Bad orbifold vector bundles have fibers given by
$\R^n / G$ for a nontrivial action of a nontrivial group $G$ over
every point of the orbifold.  Hence, as the sections of an orbifold
vector bundle are required to take values in the fixed-point
subspace of the fiber, the rank of the possible values of a section
of a bad bundle is smaller than the rank of the bundle at every
point. In some cases, such bundles admit only the zero section (see
Example \ref{ex-trivbundl} below).  Additionally, since bad orbifold
vector bundles cannot be studied as vector bundles over reduced
orbifolds, many of the natural maps defined between orbifolds cannot
be used to pull back a bad orbifold vector bundle in a well-defined
way.

In this paper, we suggest a method for extending the geometric
techniques used to study good orbifold vector bundles to the case of
bad orbifold vector bundles.  Our main result is the following
theorem.

% XXXXXXXXXX    THRM: each bad ovb is the restriction of a good ovb
% XXXXXXXXXX    THRM: each bad ovb is the restriction of a good ovb

\begin{thm}
\label{thrm-everybadarestriction}

For every bad orbifold vector bundle $\rho : E \rightarrow Q$, there
is an orbifold $R$ of which $Q$ is a suborbifold and a good orbifold
vector bundle $\hat{\rho} : \hat{E} \rightarrow R$ such that $E$ is
isomorphic to the restriction $\hat{E}_{|Q}$ of $\hat{E}$ to $Q$.

\end{thm}

As an application of this result, we demonstrate how the Chern-Weil
description of characteristic classes can be used to define
characteristic classes for bad orbifold vector bundles.  Our primary
focus is the Euler class, which we show plays the role of an
obstruction to the existence of nonvanishing sections.

In Section \ref{sec-badrtogood}, we develop the tools necessary to
prove Theorem \ref{thrm-everybadarestriction} and demonstrate an
explicit construction of a good orbifold vector bundle that
restricts to a given bad bundle $E$. Additionally, we show how this
construction allows the Chern-Weil construction of characteristic
classes.  In Section \ref{sec-applications}, we apply these results
to the Euler class and extend the obstruction property of the Euler
class to the case of bad orbifold vector bundles. The reader is
referred to \cite{satake1}, \cite{satake2}, \cite{ruangwt},
\cite{chenruan}, and \cite{mythesis} for relevant background
material.  In particular, \cite{ruangwt} contains as an appendix a
modern introduction to orbifolds;  for the most part, we follow its
notation.

Several conversations with Alexander Gorokhovsky were essential
to the development of the results contained in this paper.
As well, the author would like to thank Carla Farsi,
William Kirwin, and Kevin Manley for helpful conversations and advice.

% XXXXXXXXXXXXXXXXXXXXXXXXXXXXXXXXXXXXXXXXXXXXXXXXXXXXXXXXXXXXXX
% XXXXXXXXXX            SECTION
% XXXXXXXXXXXXXXXXXXXXXXXXXXXXXXXXXXXXXXXXXXXXXXXXXXXXXXXXXXXXXX

\section{Bad Bundles as Restrictions of Good Bundles}
\label{sec-badrtogood}

Let $\rho : E \rightarrow Q$ be a rank $k$ orbifold vector bundle
over the orbifold $Q$, which we do not assume to be reduced.  Let
$K_b$ denote the (isomorphism class of the) finite group that acts
trivially in each orbifold chart for $Q$, and let $K_f$ denote the
(isomorphism class of the) subgroup of $K_b$ that acts trivially in
orbifold vector bundle charts for $E$.  Then $E$ is bad precisely
when $K_f$ is a proper subgroup of $K_b$.

It is well-known in the case of manifolds that the the pullback via the projection
of $E$ over its own total space, $\rho^\ast E$, is isomorphic to the
vertical tangent space $VE$. For bad orbifold bundles, however, the
pullback $\rho^\ast E$ is not well-defined.  In fact, as $\Sigma_Q =
Q$, the preimage of the regular points $\rho^{-1}(Q_{reg})$ is
empty. Therefore, the projection is never a regular map, and the map
$\rho$ need not admit a unique compatible system (see \cite[Section 4.4]{ruangwt}
for more details). However, in the case of a bad orbifold vector
bundle, we do have that the vertical tangent bundle of $E$ coincides
with $E$ when restricted to the zero section.

To make this statement precise, we begin with the following definition,
similar to the manifold case.

% XXXXXXXXXX    DEFINITION: vertical vectors/bundle
% XXXXXXXXXX    DEFINITION: vertical vectors/bundle

\begin{defn}

Let $\rho : E \rightarrow Q$ be a rank $k$ orbifold vector bundle
over the orbifold $Q$.  Let $\rho_T : TE \rightarrow E$ denote the
orbifold tangent bundle of the orbifold $E$.  A vector $X \in TE$ is
{\bf vertical} if, for each lift $\tilde{X}$ of $X$ into the tangent
bundle of $V \times \R^k$ for a chart $\{ V \times \R^k, G,
\tilde{\pi} \}$ for $E$, $\tilde{X}$ is a vertical vector for the
trivial vector bundle $V \times \R^k$.  In other words,
$\tilde{X}(\rho^\ast f) = 0$ for each $G$-invariant $f \in {\mathcal
C}^\infty (V)$ (note that we use $\rho$ to denote the projection $
\rho : V \times \R^k \rightarrow V$ in a chart as well as the
projection for $E$). The collection of vertical vectors is the {\bf
vertical tangent bundle} $\rho_V : VE \rightarrow E$, where $\rho_V$
is the restriction of $\rho_T$ to $VE$.

\end{defn}

It is clear that the definition of a vertical vector does not depend
on the choice of chart, and that the set of vertical vectors forms a
sub-orbifold vector bundle of $TE$.

% XXXXXXXXXX    PROP: the restriction bndl iso to the original bndl
% XXXXXXXXXX    PROP: the restriction bndl iso to the original bndl

\begin{prop}
\label{prop-isoprop}

Let $\rho : E \rightarrow Q$ be a vector bundle over the orbifold
$Q$ as above.  Then if $\iota : Q \rightarrow E$ denotes the
embedding of $Q$ into $E$ as the zero section, the restriction
$VE_{| \iota(Q)}$ of $VE$ to $\iota(Q)$ is naturally isomorphic to
the bundle $E$.

\end{prop}

\begin{pf}

Choose a compatible cover of $Q$, and then for each injection $i :
\{ V_1, G_1, \pi_1 \} \rightarrow \{ V_2, G_2, \pi_2 \}$, there is
an associated transition map $g_i : V_1 \rightarrow AUT(\R^k)$ (see
\cite[Section 4.1]{ruangwt}). Associated to each chart $\{ V, G, \pi
\}$ in the cover is a chart $\{ V \times\R^k, G, \tilde{\pi} \}$ for
$E$, and the resulting cover is a compatible cover of $E$. Moreover,
each injection of charts for $E$ is induced by an injection in $Q$
by associating $i : \{ V_1, G_1, \pi_1 \} \rightarrow \{ V_2, G_2,
\pi_2 \}$, consisting of an embedding $\phi : V_1 \rightarrow V_2$
and an injective homomorphism $\lambda : G_1 \rightarrow G_2$, with
the injection $i_E : \{ V_1 \times \R^k, G_1, \tilde{\pi}_1 \}
\rightarrow \{ V_2 \times \R^k, G_2, \tilde{\pi}_2 \}$, consisting
of the embedding
\[
\begin{array}{lcccl}
    \phi_E  &:& V_1 \times \R^k  &\rightarrow&   V_2 \times \R^k       \\\\
            &:& (x, v)          &\mapsto&       (\phi(x), g_i(x)v)
\end{array}
\]
and homomorphism $\lambda_E := \lambda$.

Using this compatible cover, the transition maps for $TE$ are given
by the differentials of the embeddings $\phi_E$, so that for $(x, v)
\in V_1 \times \R^k$, with respect to an injection $i_E$, we have
$(g_i)_{TE} (x, v) = (d\phi_x, d(g_i(x)) ) \in AUT(\R^k) \times AUT
(\R^k) \subseteq AUT(\R^{2k})$.

The vertical vectors of $TE$ clearly correspond to elements of $V
\times T \R^k \subseteq TV \times T\R^k = T(V \times \R^k)$.  For
the bundle $VE$, then, the transition maps are $(g_i)_{VE} (x, v) =
d(g_i(x)) \in AUT(\R^k) $. Restricted to the zero section, we have
$(g_i)_{VE_{| \iota(Q)}} (x, \mathbf{0}) = d(g_i(x))$, and as each
$g_i (x)$ is linear, $d(g_i(x)) = g_i(x)$. Therefore, identifying
$(x, \mathbf{0}) \in V_1 \times \R^k$ with $x \in V_1$, we see that
the transition maps $(g_i)_{VE_{| \iota(Q)}}$ for $VE_{| \iota(Q)}$
are identical to the transition maps $g_i$ for $E$. Therefore, the
bundles are isomorphic. \qed

\end{pf}

With this, we realize the orbifold vector bundle $E$ as a
restriction of another orbifold vector bundle, $VE$, to a subset of
the base. In fact, we have the following two lemmas.

% XXXXXXXXXX    LEMMA: VE is good
% XXXXXXXXXX    LEMMA: VE is good

\begin{lem}
\label{lemm-vertgood}

The bundle $\rho_V :VE \rightarrow E$ is a good orbifold vector
bundle.  If the fiber-kernel $K_f$ for $E$ is nontrivial, then
$K_f$ acts trivially on the total space $VE$.

\end{lem}

\begin{pf}

We have that $VE$ is a subbundle of the tangent bundle $TE$ of $E$.
In each chart, the group action for the fibers of the tangent bundle
is defined to be the differential of the group action on the base.
Hence, a group element acts trivially on the fibers of $TE$ if and
only if it acts trivially on the base, and $TE$ is good.

As the kernel of the group action on the fibers of an orbifold
vector bundle is clearly contained in the kernel of the group action
on the base space, any sub-bundle of a good orbifold vector bundle
is clearly good.  Therefore, as $VE$ is a sub-bundle of $TE$, $VE$
is a good orbifold vector bundle. \qed

\end{pf}

The reader is warned that $E$ need not be a reduced orbifold.

% XXXXXXXXXX    LEMMA iota(Q) diffeo to Q
% XXXXXXXXXX    LEMMA iota(Q) diffeo to Q

\begin{lem}
\label{lemm-zerodiffeo}

The zero section $\iota(Q)$ is a suborbifold of $E$ that is
diffeomorphic to $Q$.

\end{lem}

\begin{pf}

The orbifold structure on $E$ is defined by a set of orbifold charts
$\{ V \times\R^k, G, \tilde{\pi} \}$ for $E$ induced by  orbifold charts $\{
V, G, \pi \}$ for $Q$.  Each such chart restricts to a chart $\{ V
\times \{ \mathbf{0} \} , G, \tilde{\pi}_{| V \times \{ \mathbf{0} \}} \}$
for $\iota(Q)$. It is clear that $\{ V \times \{ \mathbf{0} \} , G,
\tilde{\pi}_{| V \times \{ \mathbf{0} \}} \}$ is isomorphic to $\{ V, G,
\pi \}$ as orbifold charts and that injections of charts for $Q$
correspond bijectively to injections of charts for $\iota (Q)$.

In orbifold charts as above, $\iota$ has a well-defined ${\mathcal
C}^\infty$ lifting $\tilde{\iota} : V \rightarrow V \times \{
\mathbf{0} \}$ as the identity on $V$, taking the group
homomorphisms to be the identity. When restricted to the image
$\iota (Q)$, it is clear that this map is bijective and that its
inverse admits a well-defined ${\mathcal C}^\infty$ lifting. \qed

\end{pf}

With this, we have proven Theorem \ref{thrm-everybadarestriction}.
In the explicit construction given by Proposition
\ref{prop-isoprop}, $R = E$ and $\hat{E} = VE$.

Recall that if $U_p$ is uniformized by a {\bf chart at $p$} (i.e. a
chart $\{ V_p, G_p, \pi_p \}$ such that $G_p$ acts as a subgroup of
$O(n)$ and $p$ is the image of the origin in $V_p$) and if $q \in
U_p$, then the conjugacy class $(g)_{G_q}$ of $g \in G_q$ is said to
be equivalent to the conjugacy class $(h)_{G_p}$ of $h \in G_p$ if
there is an injection $i$ from a chart at $q$ into the chart at $p$
whose homomorphism $\lambda$ maps $g$ to $h$ (see \cite[Section
3.1]{chenruan}). Clearly, as the charts and injections for $E$ can
be taken to be those induced by its structure as a vector bundle
over $Q$, we see that the set of equivalence classes for group
elements of $Q$ coincide with those for $E$. Hence, we let $T$
denote the set of equivalence classes for both $Q$ and $E$.

The embedding $\iota : Q \rightarrow E$ induces a natural ${\mathcal
C}^\infty$ map on the spaces of sectors of these orbifolds (see
\cite[Definition 3.1.2]{chenruan}).  In particular, for each sector
$\tilde{E}_{(g)}$, let $\iota_{(g)} : \tilde{Q}_{(g)} \rightarrow
\tilde{E}_{(g)}$ be defined by setting $\iota_{(g)} [(p, (g))] =
[(p, \mathbf{0}), (g)]$.  This map has a well-defined ${\mathcal
C}^\infty$ lifting to the embedding of $V^g$ into $(V \times
\R^k)^g$ as the zero section in an orbifold chart. Hence, on each
twisted sector $\tilde{E}_{(g)}$, it induces a restriction
homomorphism
\[
    \iota_{(g)}^\ast : H_{dR}^\ast(\tilde{E}_{(g)}) \rightarrow H_{dR}^\ast(\tilde{Q}_{(g)})
\]
in the de Rham cohomology of the sector.  In the case that $Q$ and
$E$ admit almost complex structures so that $H_{orb}^\ast(Q)$ and
$H_{orb}^\ast(E)$ are defined, we let
\[
\begin{array}{rcl}
    \iota^\ast   &:=&
            \bigoplus\limits_{(g) \in T} \iota_{(g)}^\ast         \\\\
                    &:& H_{orb}^\ast(E) \rightarrow H_{orb}^\ast(Q)
\end{array}
\]
be the sum of these homomorphisms defined in Chen-Ruan orbifold
cohomology (note that throughout, we use cohomology with real coefficients).

% XXXXXXXXXX    Lemma: iota ast injective
% XXXXXXXXXX    Lemma: iota ast injective

\begin{lem}
\label{lemm-crcohomiso}

Suppose $E$ and $Q$ admit almost complex structures.  Then the
restriction homomorphism $\iota^\ast : H_{orb}^\ast(E)
\rightarrow H_{orb}^\ast(Q)$ is an isomorphism of ungraded
additive groups.

\end{lem}

The reader is warned that this isomorphism does not generally
preserve the grading; the degree shifting number of a group element
$g \in K_b$ is zero with respect to the orbifold structure on $Q$,
while it need not be zero with respect to the orbifold structure of
$E$ (see Example \ref{ex-trivbundl} below).

\begin{pf}

A compatible cover of orbifold charts for $\tilde{E}$ can be taken
as follows.  Each chart for $Q$ at a point $p$ of the form $\{ V_p,
G_p, \pi_p \}$ again induces a chart $\{ V_p \times \R^k, G_p,
\tilde{\pi}_p \}$ for $E$. Then, for each equivalence class $(g) \in
T$ with representative $g \in G_p$, there is a chart $\{ (V_p \times
\R^k)^g, C(g), \tilde{\pi}_{p, g} \}$ for a neighborhood of $[(p,
\mathbf{v}), (g)] \in \tilde{E}_{(g)}$.

Now, consider the map
\[
\begin{array}{rcccl}
    H_t     &:&     \tilde{E}_{(g)} &\rightarrow&   \tilde{E}_{(g)}  \\\\
            &:&     [(p,\mathbf{v}), (g)] &\mapsto&   [(p,t\mathbf{v}), (g)]
\end{array}
\]
for $t \in [0, 1]$.  It is clear that this map defines a ${\mathcal
C}^\infty$ deform retraction of $\tilde{E}_{(g)}$ onto the image of
$\iota_{(g)}$. Therefore, $\iota_{(g)}^\ast$ is an isomorphism
between $H_{dR}^\ast(\tilde{E}_{(g)})$ and
$H_{dR}^\ast(\tilde{Q}_{(g)})$.  As $H_{orb}^\ast (Q)$ and
$H_{orb}^\ast (E)$ are defined to be sums of the de Rham cohomology
of the sectors of the respective orbifolds, the sum, $\iota^\ast$,
is then clearly an isomorphism of additive groups. \qed

\end{pf}

Finally, we note that any characteristic class $c$ (see
\cite{ruangwt} Proposition 4.3.4) or orbifold characteristic class
$c_{orb}$ (see \cite{mythesis} Appendix B) defined by the Chern-Weil
construction for vector bundles can be extended to the case of bad
orbifold vector bundles.  Note that the definition of orbifold
characteristic classes $c_{orb}(E)$ applies the Chern-Weil
construction to connections on the bundle $\tilde{E}$ over
$\tilde{Q}$.

\begin{defn}
\label{charclass}

If $c$ is a characteristic class, we define
\[
    c(E) := \iota^\ast c(VE)
\]
to take its value in $H_{dR}^\ast(Q)$. In the case that $E$ and $Q$
admit almost complex structures, we define
\[
    c_{orb}(E)
        :=  \iota^\ast c_{orb}(VE)
        =   \iota^\ast c(\widetilde{VE})
\]
to take its value in $H_{orb}^\ast(Q)$.

\end{defn}

As these definitions are given in terms of good bundles, they
satisfy the expected properties.  The following two lemmas show that
they coincide with the original definitions when those are defined.

% XXXXXXXXXX    Lemma: char class generalizes the original definition
% XXXXXXXXXX    Lemma: char class generalizes the original definition

\begin{lem}
\label{lemm-cgb}

Let $E$ be a good orbifold vector bundle and let $c(E)$ denote a
characteristic defined by the Chern-Weil construction. Then
Definition \ref{charclass} coincides with the previous definition of
$c(E)$ : $\iota^\ast c(VE) = c(E)$.

\end{lem}

\begin{pf}

As $E$ is good, $E_{red}$ is naturally an orbifold vector bundle
over $Q_{red}$ with identical geometry to that of $E$ (see
\cite{ruangwt}). Hence, without loss of generality, we assume $Q$ is
reduced.

The projection $\rho$ is a regular map, so that it admits a unique
isomorphism class of compatible systems (see \cite[Lemma 4.4.11 and
Remark 4.4.12b]{ruangwt}). Hence, the pullback bundle $\rho^\ast E$
of $E$ is well-defined.

Note that in this compatible system, a chart $\{ V \times \R^k, G,
\tilde{\pi} \}$ for $E_{| U}$ is associated to the chart $\{ V, G,
\pi \}$ for $U \subseteq Q$ from which it was induced.  The lift
$\tilde{\rho}_{U, E_{| U}} : V \times \R^k \rightarrow V$ of $\rho$
is simply the projection onto the first factor, and the injection
$\lambda(i)$ associated to an injection
\[
    i : \{ V_1 \times\ \R^k, G_1, \tilde{\pi}_1 \} \rightarrow \{ V_2 \times \R^k,
    G_2, \tilde{\pi}_2 \}
\]
is the obvious injection
\[
    \lambda(i) : \{ V_1, G_1, \pi_1 \} \rightarrow \{ V_2, G_2,
    \pi_2 \}
\]
where the embedding is restricted to the zero section and the group
homomorphism of $\lambda(i)$ is taken to be that of $i$.

An argument similar to that used in Proposition \ref{prop-isoprop}
shows that in this case, $\rho^\ast E$ is isomorphic to $VE$ as
orbifold bundles over $E$. If $g_{\lambda(i)}$ denotes the
transition map for $E$ corresponding to an injection $\lambda(i)$
(which can be taken to be associated to an injection $i$ for $E$),
then the corresponding transition map on the vertical tangent bundle
is again given by $(g_{\lambda(i)})_{VE} (x, v) =
d(g_{\lambda(i)}(x)) \in AUT(\R^k) $. The transition maps of
$\rho^\ast E$ are the pullbacks of the transition maps via the lifts
of $\rho$ given above:
\[
    (g_i)_{\rho^\ast E}  =   g_{\lambda(i)} \circ \tilde{\rho}_{U, E_{|
    U}}.
\]
Hence, the transition maps for $\rho^\ast E$ are constant along the
fibers of $E$ and are seen to be equal to the transition maps for
$VE$.

With this, as
\begin{equation}
\label{eq-ruanlemma}
    \rho^\ast c(E) = c(\rho^\ast E)
\end{equation}
by \cite[Lemma 4.4.3]{ruangwt}, we have
\[
\begin{array}{rcl}
    \iota^\ast c(VE)    &=& \iota^\ast c(\rho^\ast E)       \\\\
                        &=& \iota^\ast \rho^\ast c(E)       \\\\
                        &&\mbox{(by Equation \ref{eq-ruanlemma})}
                                                            \\\\
                        &=& (\rho \circ \iota)^\ast c(E)    \\\\
                        &=& c(E),
\end{array}
\]
as $\rho \circ \iota$ is the identity on $Q$ up to the
diffeomorphism given by Lemma \ref{lemm-zerodiffeo}.\qed

\end{pf}

% XXXXXXXXXX    Lemma: orb char class generalizes the original definition
% XXXXXXXXXX    Lemma: orb char class generalizes the original definition

\begin{lem}
\label{lemm-corbgb}

Definition \ref{charclass} extends the original definition of
orbifold characteristic classes: if $c_{orb}$ is defined as in
\cite[Lemma 4.4.1]{mythesis} for a vector bundle $E$, then
$\iota^\ast c_{orb}(VE) = c_{orb}(E)$.

\end{lem}

\begin{pf}

In this case, for each $(g) \in T$, the bundle $ \rho_{(g)} :
\tilde{E}_{(g)} \rightarrow \tilde{Q}_{(g)}$ is a good bundle.
Hence, we can apply Lemma \ref{lemm-cgb} to see that (the reduction
of) $\tilde{E}_{(g)}$ pulls back via its projection to a bundle
isomorphic to its vertical tangent bundle.  Moreover, we have that
$c(\tilde{E}_{(g)}) = \iota^\ast c(V\tilde{E}_{(g)})$.

With this, we need only note that $c_{orb}(E)$ is the sum
$\sum\limits_{(g) \in T} c(\tilde{E}_{(g)})$, so that equality on
each sector implies that $\iota^\ast c_{orb}(VE) = c_{orb}(E)$. \qed

\end{pf}

% XXXXXXXXXXXXXXXXXXXXXXXXXXXXXXXXXXXXXXXXXXXXXXXXXXXXXXXXXXXXXX
% XXXXXXXXXX            SECTION: Applications
% XXXXXXXXXXXXXXXXXXXXXXXXXXXXXXXXXXXXXXXXXXXXXXXXXXXXXXXXXXXXXX

\section{Application: The Orbifold Euler Class of a Bad Orbifold Vector
Bundle}
\label{sec-applications}

In previous work (see \cite{mythesis} and \cite{mycyclicobstruc}),
we have defined an {\bf orbifold Euler class} in Chen-Ruan orbifold
cohomology (note that in \cite{mythesis}, the definition of an
orbifold vector bundle was taken to be that of a good orbifold
vector bundle). In \cite{mycyclicobstruc}, it was demonstrated that
when all of the local groups are cyclic, this class acts as a
complete obstruction to the existence of nonvanishing tangent vector
fields.

In this section, we show how the results of the last section can be
used to extend some of the techniques used to study good orbifold
vector bundles to the case of bad orbifold vector bundles. In
particular, we extend the definition of $e_{orb}(E)$ to the case of
$E$ bad and show that it acts as an obstruction to the existence of
nonvanishing sections.  We begin with the following lemma.

% XXXXXXXXXX    Lemma: E admits a nvs iff VE does
% XXXXXXXXXX    Lemma: E admits a nvs iff VE does

\begin{lem}
\label{lemm-nonvsec}

The bundle $\rho : E \rightarrow Q$ admits a nonvanishing section if
and only if $\rho_{VE} : VE \rightarrow E$ admits a nonvanishing
section.

\end{lem}

\begin{pf}

Suppose $VE$ admits a nonvanishing section $s$, and then the
restriction of $s$ to the zero section $\iota(Q)$ of $E$ is a
nonvanishing section of $VE_{| \iota(Q)}$.  Using the bundle
isomorphism between $VE_{| \iota(Q)}$ and $E$ given in Proposition
\ref{prop-isoprop}, $s_{| \iota(Q)}$ defines a nonvanishing section
of $E$.

Conversely, suppose $E$ admits a nonvanishing section $s$.  Define
$\hat{s}$ on $\iota(Q) \subset E$ using the isomorphism between $E$
and $VE_{| \iota Q}$ and extend $\hat{s}$ to $E$ by defining its
value to be constant along each fiber of $E$ within each chart of
the form $\{ V \times \R^k, G, \tilde{\pi} \}$.  The result is
clearly $G$-invariant, and hence defines a vertical vector field
$\hat{s}$ on $E$ that is smooth and nonvanishing. \qed

\end{pf}

With this, we extend the definition of the orbifold Euler class to
the case of bad orbifold vector bundles.

% XXXXXXXXXX    Def: orbifold euler class of bad bundle
% XXXXXXXXXX    Def: orbifold euler class of bad bundle

\begin{defn}

Let $\rho : E \rightarrow Q$ be an orbifold vector bundle, and suppose
both $Q$ and $E$ admit almost complex structures. Define
the {\bf orbifold Euler class} of $E$, $e_{orb}(E)$, to be the image
of the orbifold Euler class of $VE$ under
$\iota^\ast : H_{orb}^\ast(E) \rightarrow H_{orb}^\ast(Q)$,
\[
    e_{orb}(E) := \iota^\ast e_{orb}(VE).
\]

\end{defn}

By virtue of Lemma \ref{lemm-corbgb}, this definition extends the
original definition of $e_{orb}(E)$ in \cite{mythesis}.

We now are prepared to state and prove the  main result of this
section, that the orbifold Euler class acts as an obstruction to the
existence of nonvanishing sections of vector bundles over orbifolds.

% XXXXXXXXXX    THRM: e_{orb} an obstruction
% XXXXXXXXXX    THRM: e_{orb} an obstruction

\begin{thm}
\label{thrm-badobstruc}

Let $\rho : E \rightarrow Q$ be an orbifold vector bundle over the
orbifold $Q$.  Suppose $E$ admits a nonvanishing section $s$. Then
$e_{orb}(E) = 0 \in H_{orb}^\ast (Q)$.

\end{thm}

\begin{pf}

Suppose $\tilde{E}$ is a good bundle over $\tilde{Q}$. Let $S$
denote the trivial rank-1 subbundle of $E$ given by the linear span
of the image of $s$, and let $\tilde{S}$ denote the span of the
nonvanishing section $\tilde{s}$ of $\tilde{E}$ induced by $s$.
Endow $E$ with a metric, and then $E$ splits into $S \oplus
S^\perp$.  Choose a metric connection $\nabla$ for $E$ that respects
this product structure, and then the induced connection
$\tilde{\nabla}$ on $\tilde{E}$ respects the obvious splitting of
$\tilde{E}$ into $\tilde{S} \oplus \tilde{S}^\perp$ using the
induced metric on $\tilde{E}$.

Let $\omega$ denote the connection form for $\nabla$,
$\tilde{\omega}$ that for $\tilde{\nabla}$, and $\Omega,
\tilde{\Omega}$ the respective curvature forms (and
$\tilde{\Omega}_{(g)}$ the curvature forms restricted to the sector
corresponding to $(g) \in T$). Then with respect to a local
orthonormal frame field $(e_1, \ldots , e_k)$ for $E$ such that $e_1
\in S$, it is clear that $\omega_{i, j} = 0$ and $\Omega_{i, j} = 0$
whenever $i= 1$ or $j = 1$. Hence, $\tilde{\Omega}_{i, j} = 0$
whenever $i = 1$ or $j = 1$.

For each $(g) \in T$, the Euler curvature form that represents
$e(\tilde{E}_{(g)}) \in H_{dR}^\ast(\tilde{Q}_{(g)})$ is defined by
\[
    E(\tilde{\Omega}_{(g)}) :=
    \left\{
    \begin{array}{ll}
        \frac{1}{2^{2m} \pi^m m! }
        \sum\limits_{\tau \in S(l)} (-1)^\tau
        (\tilde{\Omega}_{(g)})_{\tau(1)\tau(2)} \wedge \cdots \wedge
        (\tilde{\Omega}_{(g)})_{\tau(l-1)\tau(l)},
            &\mbox{if} \; l
            \; \mbox{is even,}                  \\\\
        0,              &\mbox{if} \; l
            \; \mbox{is odd},
    \end{array}
    \right.
\]
where $l$ is the rank of $\tilde{E}_{(g)}$ (which may be less than
the rank $k$ of $E$ if $(g) \neq (1)$).  In this case, as at least
one of the $\tilde{\Omega}_{i,j}$ factors vanishes in each term, we
have that $E(\tilde{\Omega}_{(g)}) = 0$.  As this is true for each
$(g) \in T$, $e_{orb}(E) = 0$.  Of course, given Lemma
\ref{lemm-corbgb}, this implies that $\iota^\ast e_{orb}(VE) = 0$ as
well.

If $E$ is a bad bundle, then by Lemma \ref{lemm-nonvsec}, $VE$ admits a
nonvanishing section.  By Lemma \ref{lemm-vertgood}, $VE$ is a good
bundle, so the above argument implies that $e_{orb}(VE) = 0 \in
H_{orb}^\ast(E)$.  Hence $e_{orb}(E) = 0 \in H_{orb}^\ast(Q)$.
\qed

\end{pf}

We conclude this paper with an example illustrating Theorem
\ref{thrm-badobstruc} applied to a bundle whose ordinary Euler class
$e(E)$ is trivial yet admits only the zero section.

\begin{exmp}
\label{ex-trivbundl}

Consider the unreduced orbifold $Q$ given by a smooth manifold
diffeomorphic to $S^2$ equipped with a trivial $\Z / 3 \Z$-action.
Since every point $p \in Q$ is fixed by a nontrivial group element,
each point of the orbifold is singular.  Now, consider the orbifold
vector bundle $E$ given by equipping the trivial real rank-2 bundle
over $S^2$ with a fiber-wise $\Z / 3\Z$-action using the standard
representation of $\Z / 3\Z$ on $\R^2$ as rotations.  For the
associated orbifold vector bundle $E$, each fiber is the quotient
$\R^2/ (\Z/3\Z)$, in which the only fixed point is the zero vector.
Hence, the only section of this bundle is the zero section.

Now, consider the vertical tangent space $VE$.  It is easy to see
that the orbifold Euler class $e_{orb}(VE)$ of $VE$ as a bundle over $E$ has nonzero
components in degree $H_{orb}^{\frac{2}{3}}(E)$ and $H_{orb}^{\frac{4}{3}}(E)$ corresponding
to the twisted sectors $\tilde{E}_{(1)}$ and $\tilde{E}_{(2)}$ (which
are both diffeomorphic to $Q$, and form trivial rank-0 bundles over
$\tilde{Q}_{(1)}$ and $\tilde{Q}_{(2)}$, respectively).

In this case, the isomorphism $\iota^\ast : H_{orb}^\ast (E)
\rightarrow H_{orb}^\ast (Q)$ maps $H_{orb}^{\frac{2}{3}}(E)$ and
$H_{orb}^{\frac{4}{3}}(E)$ into $H_{orb}^0(Q)$; the degree-shifting
number is $\frac{1}{3}$ for the generator $1 \in \Z/3\Z$ with
respect to the orbifold structure of $E$, but is 0 with respect to
the orbifold structure of $Q$. Hence, the orbifold Euler class
$e_{orb}(E) := \iota^\ast e_{orb}(VE)$ of $E$ in $H_{orb}^\ast(Q)$
has nontrivial terms in degree 0.

\end{exmp}

% XXXXXXXXXXXXXXXXXXXXXXXXXXXXXX    THE BIBLIOGRAPHY
% XXXXXXXXXXXXXXXXXXXXXXXXXXXXXX    THE BIBLIOGRAPHY

\bibliographystyle{amsplain}

\end{document}